\input amstex
\documentstyle{amsppt}
\TagsOnRight \headline={\tenrm\hss\folio\hss} \magnification=\magstep1 \nopagenumbers

\topmatter
\title On the Howson property of descending HNN-extensions of
 groups \endtitle
\author Moldavanskii David I.\endauthor

\abstract A group $G$ is said to have the Howson property (or to be a Howson group) if the
intersection of any two finitely generated subgroups of $G$ is finitely generated subgroup. It
is proved that descending HNN-extension is not a Howson group under some assumptions satisfied
by the  base group of HNN-extension. In particular, a result of the paper joined with a Burns
-- Brunner result (received in 1979) implies that any descending HNN-extension of non-cyclic
free group does not have the Howson property.\endabstract
\endtopmatter
\document

\centerline{\bf 1. Introduction. Main results} \medskip

A group $G$ is said to have the Howson property (or to be a Howson group) if the intersection
of any two finitely generated subgroups of $G$ is the finitely generated subgroup too. This
denomination was introduced into practice after the work of A.~G.~Howson [1], where it was
proved that any free group possesses this property. Then, generalizing this result,
B.~Baumslag [2] have shown that a free product of two Howson groups is a Howson group. On the
other hand, it was noted in [3] that the direct product of free group of rank 2 and of
infinite cyclic group does not have the Howson property. This observation was then extended by
R.~Burns and A.~Brunner: they have proved in [4] that any extension of non-cyclic finitely
generated free group by infinite cyclic group is not a Howson group. Since every extension by
infinite cyclic group is splittable, any such group is a special case of descending
HNN-extension of free group.

Recall that descending (or named by some authors as ascending)  HNN-extension is, in
turn, a special case of general construction of HNN-extension and can be defined as
follows.

Let $G$ be a group and let $\varphi$ be an injective endomorphism of $G$. Then descending
HNN-extension of (base) group $G$ with respect to endomorphism $\varphi$ is the group
$G(\varphi)=\bigl(G, t;\ t^{-1}gt=g\varphi\ (g\in G)\bigr)$ generated by generators of $G$ and
by one more element $t$ and defined by all defining relations of $G$ and by all relations of
form $t^{-1}gt=g\varphi$ where $g\in G$.  It is obvious that if endomorphism $\varphi$ is, in
addition, surjective (i.~e. it is an automorphism of $G$) then the group $G(\varphi)$ turns
out to be a splitting extension of group $G$ by infinite cyclic group with generator $t$.
Therefore, the following assertion can be considered as a supplement to the Burns -- Brunner
above result:

\proclaim{\indent Theorem 1} Let $G$ be non-cyclic finitely generated free group and
let $\varphi$ be an injective but not surjective endomorphism of $G$. Then the
descending HNN-extension $G(\varphi)=\bigl(G, t;\ t^{-1}gt=g\varphi\ (g\in G)\bigr)$ is
not a Howson group.
\endproclaim

Thus, this result joined with the Burns -- Brunner result implies that {\it any
descending HNN-extension of non-cyclic free group does not have the Howson property}.
\smallskip

The assumption that the base group of the HNN-extension is non-cyclic is essential.
Indeed, any HNN-extension of infinite cyclic group is an one-relator group $G_k=\langle
a, t;\ t^{-1}at=a^k\rangle$ (where $k$ is a non-zero integer) belonging to the family
of Baumslag -- Solitar groups, and it was shown in [3] that all $G_k$ are Howson
groups. It is relevant to mention that this result was generalized in [4] as follows:

The group $G=\langle a_1, a_2, \dots , a_m, t;\ t^{-1}ut=v\rangle$, where $u$ and $v$
are non-identity elements of free group $F=\langle a_1, a_2, \dots , a_m\rangle$, is
Howson group provided that at least one of $u$ and $v$ is not a proper power in $F$.

One more family of one-relator Howson groups provides the result of work [5] asserting
that the generalized free product of two free groups with cyclic amalgamated subgroup
which is isolated at least in one of free factors is a Howson group.

On the other hand, many one-relator groups do not possess the Howson property. It was shown in
[3] that if non-abelian one-relator group with non-trivial center is not isomorphic to group
$G_{-1}=\langle a, t;\ t^{-1}at=a^{-1}\rangle$ then it is not a Howson group. It should be
noted that this assertion turns out to be a consequence of the Burns -- Brunner result since
non-cyclic one-relator group with non-trivial center is an extension of non-cyclic finitely
generated free group by infinite cyclic group [6]. Recently some new examples of one-relator
groups without Howson property were given in [7], [8] and [9]. However, it is easy to see that
all these groups are a descending HNN-extensions of non-cyclic free group. Thus, the
impracticability of Howson property in all examples of one-relator non-Howson groups that we
know up to now is in fact a consequence of our Theorem 1 and Burns -- Brunner
result.\smallskip

Theorem 1 is a special case of the following somewhat more general result. Let us say that a
subgroup $H$ of group $G$ is freely complemented if  there exists a non-identity subgroup $K$
of $G$ such that subgroup generated by subgroups $H$ and $K$ is their free product $H\ast K$.

\proclaim{\indent Theorem 2} {\it Let $G$ be a finitely generated group, let $\varphi$ be an
injective but not surjective endomorphism of $G$ and $H=G\varphi$. If subgroup $H$ of group
$G$ is freely complemented then the descending HNN-extension $G(\varphi)$ is not a Howson
group.
\endproclaim

In order to deduce Theorem 1 from Theorem 2 it is enough to note that if $G$ is a
non-cyclic finitely generated free group and $\varphi$ is an injective but not
surjective endomorphism of $G$ then subgroup  $H=G\varphi$ is freely complemented. In
fact, since rank of subgroup $H$ is equal to rank of $G$ and $H$ is a proper subgroup
of $G$, the Schreier's formula implies that $H$ is of infinite index in $G$. Therefore,
it follows from the Hall -- Burns Theorem (see e.~g. [10, proposition 1.3.10]) that $H$
is freely complemented.

One more application of Theorem 2 is

\proclaim{\indent Corollary} Let a finitely generated group $G$ is the free product of
non-identity groups $A$ ¨ $B$. If $\varphi$ is an injective but not surjective endomorphism of
$G$ such that $A\varphi\subseteq A$ and $B\varphi\subseteq B$ then $G(\varphi)$ is not a
Howson group.
\endproclaim

In this case subgroup $H=G\varphi$ is generated by subgroups $A\varphi$ ¨ $B\varphi$
(and is a free product of them) and since $H\varphi\neq G$ then $A\varphi\neq A$ or
$B\varphi\neq  B$. Therefore, $H$ is of infinite index in $G$ and hence (see e.~g. [11,
p.~27]) subgroup $H$ is freely complemented.\medskip

The similar assertion is fulfilled for group that is decomposable into the direct
product:

\proclaim{\indent Theorem 3} Let group $G$ be a direct product of non-identity groups
$A$ and $B$ and let $\varphi$ be an injective but not surjective endomorphism of $G$
such that $A\varphi\subseteq A$ ¨ $B\varphi\subseteq B$. If $A\varphi\neq A$,
$B\varphi\neq B$ and at least one of subgroups $A$ and $B$ is finitely generated then
$G(\varphi)$ is not a Howson group.
\endproclaim

Theorem 3 implies, in particular, that if $G$ is a free abelian finitely generated group and
$\varphi$ is injective endomorphism of $G$ such that the matrix of $\varphi$ in some free base
of $G$ is of block-diagonal form where determinant of at least two diagonal blocks is not
equal to $\pm 1$ then $G(\varphi)$ is not a Howson group. The problem of complete
characterization of those descending HNN-extensions of free abelian groups that are a Howson
groups is still open.\medskip

\centerline{\bf 2. The proof of Theorem 2}
\medskip

Let $\varphi$ be an injective but not surjective endomorphism of finitely generated group $G$,
let $H=G\varphi$ and $K$ be an non-identity subgroup of $G$ such that subgroup $L$ generated
by subgroups $H$ and $K$ is their free product, $L=H\ast K$. It is obvious that we can assume
subgroup $K$ to be finitely generated.

For any integer $n$ let $K_{n}=t^{-n}Kt^{n}$. Let also $N$ denote the subgroup of group
$G(\varphi)$ that is generated by all subgroups $K_{n}$ and $M$ denote the subgroup of group
$G(\varphi)$ that is generated by all subgroups $K_{n}$ with $n\geqslant 0$. Remark that for
$n\geqslant 0$ we have $K_{n}=K\varphi^{n}$ and therefore subgroup $M$ is contained in the
base group $G$ of  HNN--extension $G(\varphi)$.

\proclaim{\indent Lemma 1} Subgroup $N$ is the free product of family subgroups $K_{n}$,
$n\in\Bbb Z$. Hence subgroups $N$ and $M$ are not finitely generated.
\endproclaim

In order to prove Lemma 1 it is enough to prove that any subgroup generated by a finite family
of subgroups $K_{n}$ is the free product of these subgroup, and to this end, in turn, it is
enough to prove that  for any integer $r\geqslant 1$ subgroup $M_r$ generated by subgroups
$K_{0}=K\varphi^0$, $K_{1}=K\varphi$, \dots , $K_{r}=K\varphi^r$ is the free product of these
subgroups.

When $r=1$ this is  obvious since $L=H\ast K=H\ast K_{0}$ and $K\varphi\leqslant H$. Let us
assume that for some $r\geqslant 1$ subgroup $M_r$ is the free product of subgroups $K_{0}$,
$K_{1}$, \dots , $K_{r}$. Then since the mapping $\varphi$ is an isomorphism of group $G$ on
the group $H$ and for any $i\geqslant 0$ $K_{i}\varphi=K_{i+1}$ subgroup $M_r\varphi$ is the
free product of subgroups $K_{1}$, $K_{2}$, \dots , $K_{r+1}$. Since subgroup $M_{r+1}$ is
generated by subgroups $K_{0}$ and $M_r\varphi$ and $M_r\varphi\leqslant H$ this implies that
subgroup $M_{r+1}$ is the free product of subgroups $K_{0}$, $K_{1}$, \dots , $K_{r+1}$. The
proof of Lemma 1 is complete.

\proclaim{\indent Lemma 2} $N\cap G=M$.
\endproclaim

 Since the inclusion $M\subseteq N\cap G$ is trivial it is enough to
prove the opposite inclusion. Any non-identity element $u$ of subgroup $N$ can be written in
the form
$$
u=v_1v_2\cdots v_r,
$$
where $r\geqslant 1$, for any $i=1, 2, \dots , r$ $v_i$ is non-identity element from some
subgroup $K_{n_i}$, $v_i=t^{-n_i}g_it^{n_i}$ for some non-identity element $g_i\in K$, and if
$r>1$ then for any $i=1, 2, \dots , r-1$ $n_i\neq n_{i+1}$.

We shall show that if at least one of the numbers $n_{1}$, $n_{2}$, \dots , $n_{r}$ is
negative, then element $u$ does not  enter in subgroup $G$. Since otherwise the inclusion
$u\in M$ is evident by that the proof of Lemma will be complete.

So, let us suppose that for some $i$, $1\leqslant i\leqslant r$, we have $n_{i}<0$. If $r=1$
then since element $g_{1}$ does not belong to subgroup $H$, the expression
$u=t^{-n_{1}}g_{1}t^{n_{1}}$ is reduced in $HNN$-extension $G(\varphi)$ and therefore $u\notin
G$ by Britton's Lemma.

Now, let $r>1$ and $n$ denote the smallest from integers $n_1$, $n_2$, \dots , $n_r$. Suppose
by the contrary that element $u$ belongs to subgroup $G$. Then since $n\leqslant -1$ element
$t^{n}ut^{-n}=u\varphi^{-n}$ belongs to subgroup $H$.

On the other hand since $n-n_i\leqslant 0$ for any $i=1, 2, \dots , r$, we have for every such
number $i$
$$
t^{n}v_it^{-n}=t^{n-n_i}g_it^{-(n-n_i)}=g_i\varphi^{n_i-n}\in K\varphi^{n_i-n}.
$$
Therefore, since for any $i=1, 2, \dots , r-1$ integers $n_i-n$ ¨ $n_{i+1}-n$ are different,
the following expression  of element $t^{n}ut^{-n}$,
$$
t^{-n}ut^n=g_1\varphi^{n_1-n}\cdot g_2\varphi^{n_2-n}\cdot \cdots\cdot g_r\varphi^{n_r-n},
$$
is reduced in decomposition of group $M$ into free product in Lemma 1.

By the choice of integer $n$ there exists at least one number $i$ such that $n_i-n=0$; let
$i_{1}<i_{2}< \cdots <i_{s}$ be all numbers of those syllables $g_i\varphi^{n_i-n}$ for which
this  equality is satisfied. The rest syllables in this expression  of element $t^{n}ut^{-n}$
belong to subgroup $H$ and by join all such consecutive syllables we obtain the expression of
element $t^{n}ut^{-n}$ of form
$$
t^{n}ut^{-n}=w_{0}g_{i_{1}}w_{1}g_{i_{2}}w_{2}\dots w_{s-1}g_{i_{s}}w_{s},
$$
where all $w_{j}$ are elements of subgroup $H$ that are not  equal to identity except for, may
be, $w_{0}$ ¨ $w_{s}$. In any case this expression is reduced in free decomposition  $L=H*K$
of subgroup $L$ and since at least one syllable of it belongs to subgroup  $K$, this
contradicts to inclusion $t^{n}ut^{-n}\in H$. Lemma 2 is proved.\smallskip

Now we can complete the proof of Theorem 2. Let $F$ be subgroup of group $G(\varphi)$
generated by subgroup $K$ and element $t$. We shall show that $F\cap G=M$. Since subgroups $F$
and $G$ are finitely generated while subgroup $M$ (by Lemma 1) is not finitely generated, this
will imply that the group $G(\varphi)$ is not a Howson group.

Arbitrary element $f\in F$ can be written in the form
$$f=g_0t^{n_1}g_1t^{n_2}\cdots t^{n_r}g_r$$ where $g_0$, $g_1$, \dots , $g_r$ are some elements
from subgroup $K$ and  $n_1$, $n_2$, \dots, $n_r$ are some integers. The factorization of
group $G(\varphi)$ by the normal closure of subgroup $G$ shoes evidently that if element $f$
belongs to subgroup $G$ then $n_1+n_2+ \cdots +n_r=0$ and therefore  $f\in N$. Thus, we have
inclusion $F\cap G\subseteq N$ and this with taking into account of Lemma 2 and obvious
inclusion $M\subseteq F$ implies that
$$
F\cap G=F\cap G\cap N=F\cap M=M.
$$

\bigskip\newpage

\centerline{\bf 3. The proof of Theorem 3} \medskip

Let $G=A\times B$ and let $\varphi$ be an injective endomorphism of group $G$ such that
$A\varphi\subseteq A$ and $B\varphi\subseteq B$. Suppose also that $A\varphi\neq A$,
$B\varphi\neq B$ and subgroup $A$ is finitely generated. Let the restriction of mapping
$\varphi$ on subgroup $B$ be denoted by $\varphi$ too and let $B(\varphi)=\bigl(B, t;\
t^{-1}bt=b\varphi\ (b\in B)\bigr)$ be corresponding descending $HNN$-extension of group
$B$.

It is easy to see that there exists a homomorphism $\rho$ of group $G(\varphi)$ to the group
$B(\varphi)$ which sends the stable letter of group $G(\varphi)$  onto stable letter of group
$B(\varphi)$ and action of which on subgroup $G$ coincides with action of projection $\pi:G\to
B$. We claim that the kernel of $\rho$ is equal to subgroup
$U=\bigcup_{k=0}^{\infty}t^{k}At^{-k}$.

Indeed, since $A\rho=A\pi=1$ the inclusion $U\subseteq \text{Ker}\,\rho$ is evident.
Backwards, arbitrary element $v$ from $\text{Ker}\,\rho$ (just as any element of group
$G(\varphi)$) can be written in form $v=t^{m}gt^{-n}$ for some integers $m\geqslant 0$ and
$n\geqslant 0$ and some element $g\in G$. Let $g=ab$ where $a\in A$ and $b\in B$. Then
$u\rho=t^{m}bt^{-n}$ and therefore in group $B(\varphi)$ we have the equality
$t^{m}bt^{-n}=1$. Since in any HNN-extension the stable letter generates subgroup that
intersects the base group trivially then  $b=1$ and $m=n$. Thus, $v=t^{m}at^{-m}\in U$ and the
proof of equality $\text{Ker}\,\rho = U$ is complete.

Remark that since $A\varphi\neq A$ subgroup $U$ is the union of strictly increasing sequence
of subgroups and therefore is not finitely generated.

Let $C$ denote subgroup of group $G(\varphi)$ generated by subgroup $A$ and element $t$ and
let $D$ denote subgroup generated by subgroup $A$ and element $tb$ where $b\in B\setminus
B\varphi$.

It is evident that $U\leqslant C$ and it is easy to see that subgroup $U$  also is contained
in $D$. In fact, we have $A\leqslant D$. If for some $k\geqslant 0$ subgroup $t^{k}At^{-k}$ is
 contained in $D$ then $D$ contains subgroup $(tb)\,t^{k}At^{-k}(tb)^{-1}$. But since
$bt^{k}=t^{k}\,b\varphi^{k}$ we have $(tb)\,t^{k}At^{-k}(tb)^{-1}=t^{k+1}At^{-(k+1)}$.

Thus, subgroup $U$ is contained in intersection of subgroups $C$ ¨ $D$. We shall prove now
that, in fact,  $C\cap D=U$. Since subgroups $C$ and $D$ are finitely generated and subgroup
$U$ is not finitely generated then the proof of Theorem 3 will be complete.

The image of subgroup $C$ under homomorphism $\rho$ of group $G(\varphi)$ on group
$B(\varphi)$ is the cyclic subgroup generated by element $t$ and the image of subgroup $D$ is
the cyclic subgroup generated by element $tb$. If the intersection $C\rho\cap D\rho$ of these
subgroups would be non-trivial then for some non-zero integers $m$ and $n$ in group
$B(\varphi)$ must be fulfilled the equation $t^{m}=(tb)^{n}$. The passage to the quotient of
group $B(\varphi)$ by normal closure of subgroup $B$ shoes that $m=n$. Consequently, in group
$B(\varphi)$ the equation $t^{m}=(tb)^{m}$ is fulfilled, where the integer $m$ may be supposed
to be positive. Since
$$
(tb)^{m}=t^{m}\cdot t^{-(m-1)}bt^{m-1}\,t^{-(m-2)}bt^{m-2}\,\cdots\, t^{-1}bt\, b=t^{m}\cdot
b\varphi^{m-1}\, b\varphi^{m-2}\,\cdots\, b\varphi\, b,
$$
we have the equality $b\varphi^{m-1}\, b\varphi^{m-2}\,\cdots\,b\varphi\, b=1$. This implies
the inclusion $b\in B\varphi$ which contradicts to the choice of element $b$.

So, $C\rho\cap D\rho=1$ and therefore $(C\cap D)\rho=1$. Since the kernel of $\rho$ coincides
with subgroup $U$ this implies the required inclusion $C\cap D\subseteq U$.

\bigskip

\centerline{\bf References}
\medskip

\noindent\item{1.} Howson A.~G., On the intersection of finitely generated free groups, J.
London Math. Soc. 29. 1954. 428--434.

\noindent\item{2.} Baumslag B., Intersections of finitely generated subgroups in free
products,  J. London Math. Soc. 41, 1966, 673--679.

\noindent\item{3.} Moldavanskii D.~I., On the intersection of finitely generated subgroups,
Siberian Math. J. 1968, ü 9, 1422--1426 (Russian).

\noindent\item{4.} Burns R., Brunner A., Two remarks on group Howson property, Algebra and
Logic, 18, ü 5 (1979), 513--522 (Russian).

\noindent\item{5.} Burns R., On finitely generated subgroups of an amalgamated product of two
subgroups, Trans. Amer. Math. Soc. 169 (1972), 293--306.

\noindent\item{6.} Baumslag G., Taylor T., The center of groups with one defining relator,
Math\. Ann. 1968, Vol.~175, 315--319.

\noindent\item{7.} Kapovich I., Howson property and one-relator groups, Communs in algebra,
1999, Vol. 27, 1057--1072.

\noindent\item{8.} Bezverhnyaya N.~B., On the Howson property and hyperbolicity of some
2-generated 1-related groups, Algorithmic problems in groups and semigroups, Tula,
 2001 (Russian).

\noindent\item{9.} Bezverhnyaya N.~B. On the Howson property of some one-relator groups,
 Chebyshevskii sbornik (Publishers of Tula State Pedagogic
University), Vol.~2 (2001), 14--18 (Russian).

\noindent\item{10.}  Lyndon R.~C.,  Schupp P.~E., Combinatorial group theory,
Springer~-~Verlag, Berlin, etc. 1977.

\noindent\item{11.} Goryushkin A.~P., Goryushkin V.~A. Free decomposability and  balance of
groups, Publishers of Kamchatka State Pedagogic University, 2005.(Russian)\bigskip

Ivanovo State University \smallskip

{\it E-mail address}: moldav\@mail.ru
\end